\documentclass[12 pt]{article}
\usepackage{amssymb,amsfonts,amsthm}
\usepackage{amsmath}
\usepackage{graphicx}
\usepackage{epsfig,subfigure}
\usepackage{bm, float}
\floatstyle{plain}

\usepackage{amssymb}

\usepackage{graphicx}
\usepackage{float}

\usepackage{bm}
\usepackage[ruled]{algorithm2e}

\usepackage[utf8]{inputenc}
\usepackage[T1]{fontenc}
\usepackage[edges]{forest}
\usetikzlibrary{fit}

\theoremstyle{definition}

\theoremstyle{remark}

\numberwithin{equation}{section}

\begin{document}

\title{A Finite Difference Method on Irregular Grids with Local Second Order Ghost Point Extension for Solving Maxwell's Equations Around Curved PEC Objects}

\author{Haiyu Zou\footnotemark[1] \; and  Yingjie Liu\footnotemark[2]}

\date{}

\maketitle

\renewcommand{\thefootnote}{\fnsymbol{footnote}}
\footnotetext[1]{({\tt E-mail: zou@gatech.edu})
\\School of Mathematics, Georgia Institute of Technology,
Atlanta, GA 30332. Research supported in part by NSF grant DMS-1522585, DMS-CDS\&E-MSS-1622453.}

\footnotetext[2]{({\tt E-mail: yingjie@math.gatech.edu})
\\School of Mathematics, Georgia Institute of Technology,
Atlanta, GA 30332. Research supported in part by NSF grant DMS-1522585, DMS-CDS\&E-MSS-1622453.}
\renewcommand{\thefootnote}{\arabic{footnote}}

\begin{abstract}
A new finite difference method on irregular, locally perturbed rectangular grids has been developed for solving 
electromagnetic waves around curved perfect electric conductors (PEC). This method incorporates the back and forth error compensation and
correction method (BFECC) and level set method to achieve convenience and higher order of accuracy at complicated PEC boundaries. 
A PDE-based local second order ghost cell extension technique is developed based on the level set framework in order to compute the 
boundary value to first order accuracy (cumulatively), and then BFECC is applied to further improve the accuracy 
while increasing the CFL number. Numerical experiments are conducted to validate the properties of the method.
\end{abstract}


\section{Introduction}
Despite the plethora of strategies devised to solve electromagnetic problems over the last several decades, the finite-difference time-domain (FDTD) method (or Yee scheme \cite{yeeNumericalSolutionInital1966}) remains to be one of the most widely used due to its simplicity and efficiency. FDTD works the best when interfaces under study are aligned well with orthogonal grids \cite{yeeNumericalSolutionInital1966, liuOverlappingYeeFDTD2009}. However, when it comes to modeling surfaces of complex objects, to strike a balance between model complexity and accuracy/stability is no easy task. Staircasing \cite{umashankarCalculationExperimentalValidation1987a} as an approximation method works well in certain modelling scenarios but suffers from stability issues and, in extreme cases, accuracy fallout \cite{liuOverlappingYeeFDTD2009}. A local conformity method proposed by Fang J. and Ren J. \cite{fangLocallyConformedFinitedifference1993} analyzes areas where curved surfaces intersect with orthogonal grids. And by using nonorthogonal grids, Jin-fa-Lee et al. \cite{jin-fa-leeModelingThreedimensionalDiscontinuities1992} have reformulated the FDTD approach by introducing two more 
unknowns in setting up its global coordinate system and obtained good body-fitting results. 
The overlapping Yee FDTD method \cite{liuOverlappingYeeFDTD2009} uses dual overlapping and body-fitting grid to achieve second order accuracy and robustness in the computation of electromagnetic waves around complicated dielectric objects.
Other approaches that are based on or inspired by finite volume methods exploit the integral formulation 
of the problem \cite{thngEdgeelementTimedomainMethod1994, ruey-beeiwuHybridFinitedifferenceTimedomain1997} and usually involve higher cost in 
implementation and computation than a finite difference scheme does. 
An important boundary condition in the study of Maxwell's equations is PEC, or Perfect Electric Conductor and its magnetic counterpart PMC - Perfect Electric Conductor \cite{seniorApproximateBoundaryConditions1995}. The PEC condition can be described as 

\begin{equation*}
	\bm{E} \times \bm{n} = 0 \;\; {\rm and} \;\; \bm{H} \cdot \bm{n} = 0,
\end{equation*}

where $\bm{n}$ is the local normal direction of the interface. Recent formulations have managed to incorporate them into a more general framework of analysis \cite{lindellElectromagneticWaveReflection2017} and they continue to be central in the theoretical studies and application of disciplines such as antenna designs.
The use of image theory in setting up an equivalent problem as a treatment of the PEC condition is well discussed in literature \cite{bannisterImageTheoryElectromagnetic1982} \cite{balanisAntennaTheoryAnalysis2016}, which mainly focuses on radiation problems over a half plane as the PEC interface. For curved surfaces, there are recent works that extend the same understanding to spheres \cite{lindellElectromagnetostaticImageTheory2006} or even more general geometry \cite{prudencioGeneralizedImageMethod2013}. 
It is used to capture PEC/PMC boundary conditions in constructing FDTD based schemes in \cite{tayImplementationsPMCPEC2010}. A finite element method for solving Maxwell's equations around PEC corners and highly curved surfaces can be found in \cite{Boyse97}.

This paper proposes a new method to solve Maxwell's equations that involve objects with potentially complicated geometries coupled with perfect electric conductor (PEC) boundary conditions. At its core, it is driven by a simple first order finite difference scheme with the boundary conditions handled by a new higher order ghost value
extension technique based on the level-set method \cite{OsherSethian88} and related redistancing \cite{sussmanLevelSetApproach1994} and extension \cite{AS98, Peng99, Fedkiw99} methods as well as a locally conforming point-shifted method \cite{McBryan80}. Then the accuracy and stability in the interior of the computational domain as well as along the PEC 
boundary are further improved by the back and forth error compensation and correction (BFECC) method \cite{dupontBackForthError2003, dupontliu07, wangBackForthError2019}.
The level set method provides a convenient platform for capturing complex interfaces \cite{sussmanLevelSetApproach1994} as well as extrapolating certain quantities
across an interface using PDE-based methods.
BFECC can be used to solve a hyperbolic system and improve order of accuracy of an underlying scheme if it is odd-order accurate. It can also improve the CFL number of
an underlying scheme or even stabilize an unstable underlying scheme \cite{dupontliu07, wangBackForthError2019}. These properties are very helpful because if the
PEC boundary is complex, locally second order accurate ghost values outside the PEC boundary are already nontrivial to obtained as will be seen later. 
This results in first order accuracy (cumulatively) in computed boundary values,  which can then be improved to second order accuracy by BFECC.

\section{Method}

	We will first cover the necessary background knowledge to implement our method.  

\subsection{Back and Forth Error Compensation and Correction for Accuracy and Stability}
\label{Sec:BFECC}
	
Back and Forth Error Compensation and Correction, or BFECC, can be used to solve homogeneous linear hyperbolic PDE systems with constant coefficients with improved accuracy and stability. To illustrate, consider, on a uniform orthogonal grid, a linear hyperbolic PDE system of the form:

	\begin{align}\label{eq:hyperbolic_system_forward}
	\partial_t \bm{u} + \sum_{i = 1}^{d} A_i \partial_{x_i} \bm{u} = 0,
	\end{align}

	and its time-reversed system

	\begin{align}\label{eq:hyperbolic_system_backward}
	\partial_t \bm{u} - \sum_{i = 1}^{d} A_i \partial_{x_i} \bm{u} = 0,
	\end{align}

	 where $A_i$ is a real constant matrix. Let $U^n$ be the numerical solution to the system at time $t_n$ and $\mathcal{L}$ be a linear numerical scheme that evolves the numerical solution from $t_n$ to $t_{n+1}$. That is,

	\begin{equation*}
	\bm{U}^{n+1} = \mathcal{L} \bm{U}^{n}
	\end{equation*}

	and let $\mathcal{L}^*$ be the operator obtained by applying the same scheme to the time reversed equation (\ref{eq:hyperbolic_system_backward}) the same time step size, thus 
	
	\begin{equation*}
	\bm{U}^{n} \approx \mathcal{L}^* \bm{U}^{n+1}.
	\end{equation*}

	Then the general method of BFECC can be summed up in these following steps. 

	\begin{enumerate}
	\item{\bf Solve forward}. \\
	$\tilde{\bm{U}}^{n+1} = \mathcal{L} \bm{U}^{n}$.
	\item{\bf Solve backward}. \\
	$\tilde{\bm{U}}^{n} = \mathcal{L}^{*} \tilde{\bm{U}}^{n+1}$.
	\item{\bf Solve forward with the modified solution at time $t_n$}. \\
	$\bm{U}^{n+1} = \mathcal{L} \left( \bm{U}^n + \bm{e} \right)$, where $\bm{e} = \frac{1}{2} \left( \bm{U}^n - \tilde{\bm{U}}^n \right)$.
	\end{enumerate}

	where $\tilde{\bm{U}}^{n}$ is an intermediate numerical solution at time $t_n$. Observe how the error is captured by $\bm{e}$ and is later used to compensate $\bm{U}^{n}$ as an initial state to improve the accuracy of the overall scheme. Generally speaking, BFECC can improve not only the order of accuracy by one for odd order underlying schemes, but also their stabilities in the sense that if the amplification factor of the underlying scheme is no more than 2, then it becomes stable after applying BFECC \cite{wangBackForthError2019}.

	\subsection{Local Approximation with the Least Square Method}

	Consider the Maxwell equations in two dimensions:

	\begin{align}\label{eqn:2d_maxwell_tmz}
	\begin{split}
	& \frac{\partial H_x}{\partial t} = - \frac{\partial E_z}{\partial y} \\
	& \frac{\partial H_y}{\partial t} =  \frac{\partial E_z}{\partial x} \\
	& \frac{\partial E_z}{\partial t} = \frac{\partial H_y}{\partial x} - \frac{\partial H_x}{\partial y}.
	\end{split}
	\end{align}

	Following \cite{wangBackForthError2019}, to solve the system on irregular grids as a result of point shifting, we resort to the least square linear fitting method. To illustrate, consider approximating $E_z$ at grid point $(x^i, y^j)$. Since the topological feature of the point-shifted grids should be mostly intact, we can still identify the neighboring points of $(x^i, y^j)$, given as $(x^{i+1}, y^{j}), (x^{i-1}, y^{j}), (x^{i}, y^{j+1}) \text{ and } (x^{i}, y^{j+1})$. These 5 points including the central one $(x^i, y^j)$ form the basis of our local approximation. Next, we would like to obtain constants $c_0, c_1 \text{ and } c_2$  such that surface $E_z(x,y) = c_0(x - x^i) + c_1(y-y^i) + c_2$ encompasses all 5 points. This is an overfitting problem where the least square method comes in handy. We can therefore approximate $E_z$ at grid point $(x^i, y^j)$ which can be used further to approximate spatial derivatives of $E_z$ by considering the gradient of the fitted linear polynomial with respect to $x$ or $y$.

	With this in mind, we discretize the system as follows.

	\begin{align}\label{scheme:ls_cd}
	\begin{split}
	\left( H_x \right)^{n+1}_{i, j} = & \left( \tilde{H}_x \right)^{n}_{i, j} - \Delta t \left( \frac{\partial \tilde{E}_z}{\partial y} \right)_{i, j}^n  \\
	\left( H_y \right)^{n+1}_{i, j} = & \left( \tilde{H}_y \right)^{n}_{i, j} + \Delta t \left( \frac{\partial \tilde{E}_z}{\partial x} \right)_{i, j}^n\\
	\left( E_z \right)^{n+1}_{i, j} = & \left( \tilde{E}_z \right)^{n}_{i, j} + \Delta t \left(  \left( \frac{\partial \tilde{H}_y}{\partial x} \right)_{i, j}^n - \left( \frac{\partial \tilde{H}_x}{\partial y} \right)_{i, j}^n \right)  \\
	\end{split}
	\end{align}

	where terms with tildes such as $\frac{ \partial \tilde{E}_z}{\partial y}$ are least squared approximations of the spatial derivatives. One can notice that this discretization is essentially a hybrid scheme based on a central difference scheme and the Lax–Friedrichs Method \cite{wangBackForthError2019}.

	We also use the above linear least square approximation to obtain the local normal direction $\bm{n}$ from the signed distance function $\phi$ to be covered in the next section.

\subsection{Level Set Method and Point Shift Method to Capture the PEC Boundary Condition}

	There are essentially two challenges in capturing a PEC boundary condition of an arbitrary geometric shape: the mismatch between the geometric shape and orthogonal grid points and the PEC condition itself. In our approach, the rectangular grid is locally perturbed by the point-shifting method \cite{McBryan80} to match the PEC boundary and a local second order ghost point extension technique is developed based on the level set method \cite{OsherSethian88, sussmanLevelSetApproach1994, AS98, Peng99, Fedkiw99} for approximating the boundary conditions. 

\begin{figure}[H]
	\centering
	\parbox{6cm}{
	\includegraphics[width=6cm]{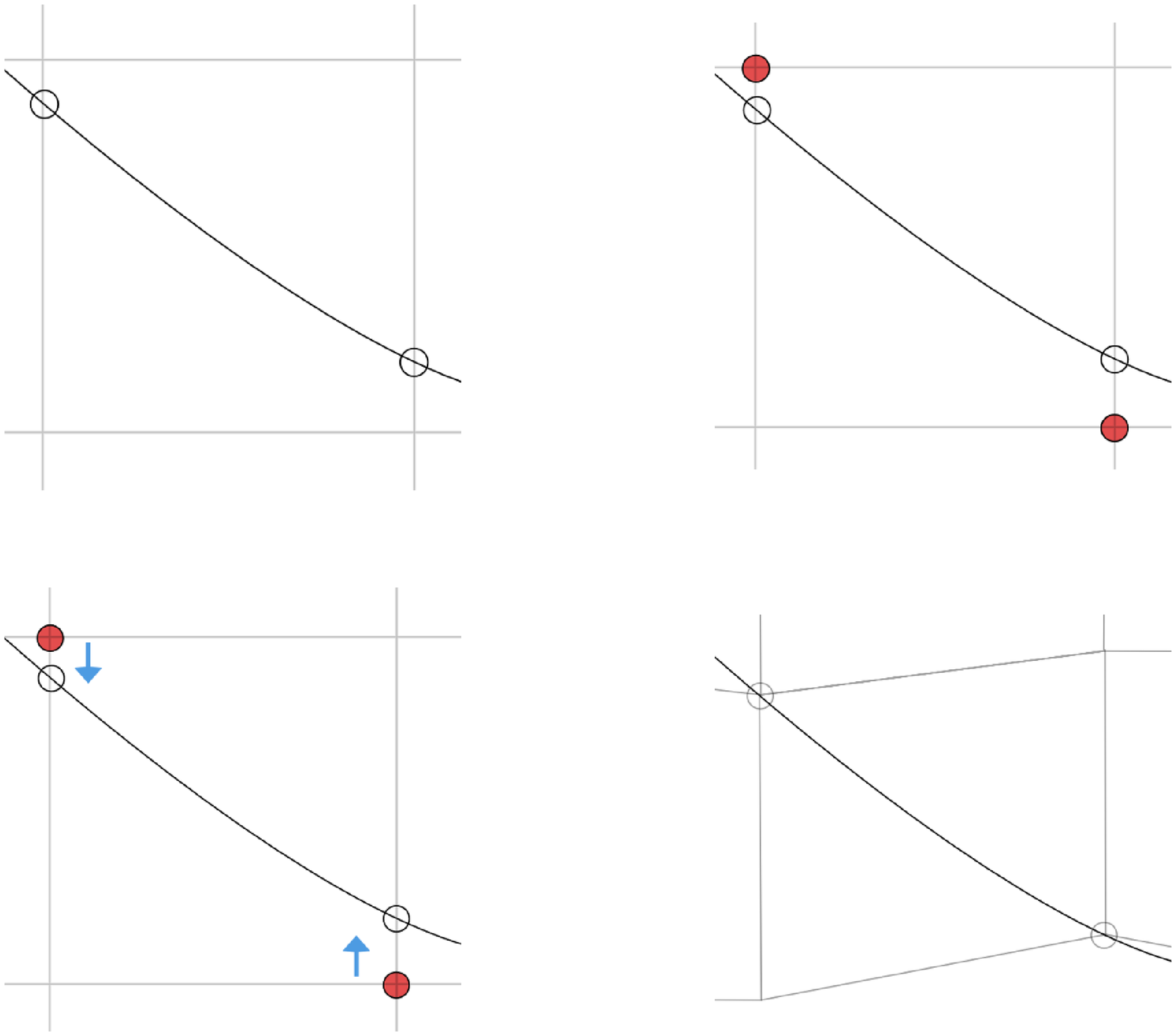}
	\caption{Constructing a conforming grid by shifting points}
	\label{fig:Point-shift-demo}}
	\qquad\qquad
	\begin{minipage}{6cm}
	\includegraphics[width=6cm]{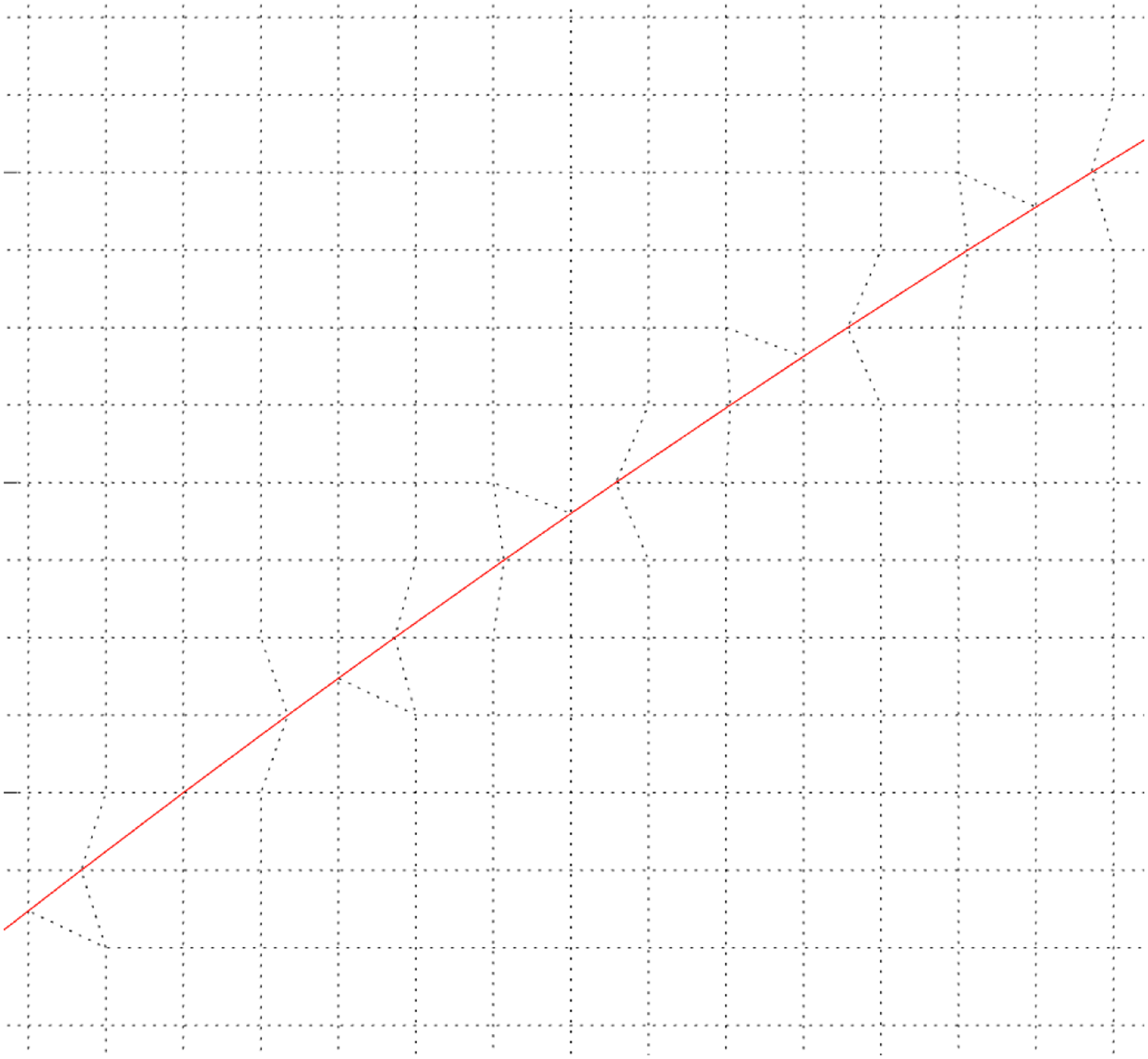}
	\caption{Grid points shifted near the interface}
	\label{fig:Point-shift}
	\end{minipage}
	
	\end{figure}

	Consider a geometric shape denoted by a closed curve $P$ in 2D with an interior, superimposed over a uniform orthogonal grid $G = H \times V \text{ where } H = \{x_i: x_i = i\Delta x, i = 1, \ldots , N_x \} \text{ and } V = \{y_j: y_i = j\Delta y, j = 1, \ldots, N_y\}$. Denote by $C$ the collection of $P$'s intersection points with the grid's lattice. That is $C = \{(x,y) \in P: x \in H \text{ or } y \in V\}$. We then go through every point in $C$, identify the closest point in $G$ and replace the grid point in $G$ with the coordinates of its corresponding intersection point. The updated grid $\tilde G$ therefore conforms to the shape of $P$. Figure \ref{fig:Point-shift-demo} gives a visual representation of this process. Figure \ref{fig:Point-shift} shows how the orthogonal grid conforms to the shape of an interface, shown in red, after point-shifting.

\begin{figure}[H]
	\centering
	\includegraphics[width=0.8\linewidth]{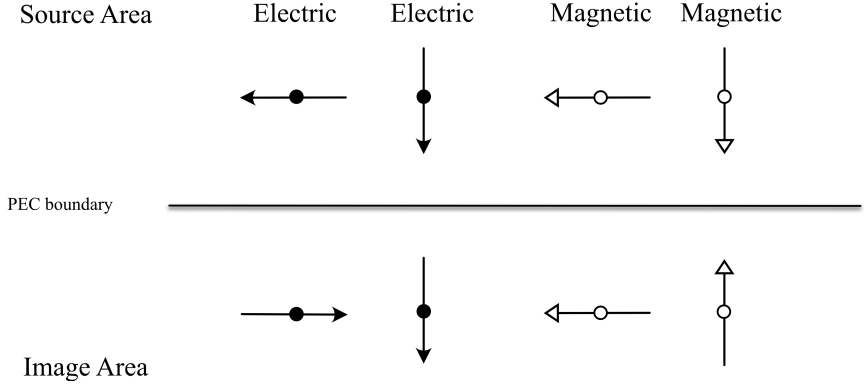}
	\caption{Mirrored fields by image theory \cite{balanisAntennaTheoryAnalysis2016}}
	\label{fig:ghost fields by image theory}
	\end{figure}

	To enforce PEC boundary condition, we use image theory to help us simplify the construction. By image theory, we set up image fields across the interface that mimics the fields under consideration \cite{balanisAntennaTheoryAnalysis2016}. See Figure \ref{fig:ghost fields by image theory}.

	However, what if the PEC interface is not well aligned with an orthogonal grid, or worse yet that it itself is not a straight line? We would like still to construct ghost fields with the help of the level set method. 

The signed distance function $\phi$ w.r.t. the PEC boundary (with positive values inside the PEC object) can be obtained by evolving the following equation until equilibrium \cite{sussmanLevelSetApproach1994},

	\begin{align}\label{eqn:distance_level_set}
	\begin{split}
	\frac{\partial \phi}{\partial t} + \bm{v} \cdot \nabla \phi = sgn (\phi)~,\\
	\end{split}
	\end{align}

	where $\bm{v} = sgn(\phi) \frac{\nabla \phi}{\| \nabla \phi \|}$. To smooth out the discontinuity over the interface, we use signed function $sgn(x)$ defined by

	\begin{align}\label{smoothed-sign-funct}
	\begin{split}
	sgn (x) = \frac{x}{\sqrt{x^2 + \Delta x^2}}\\
	\end{split}
	\end{align}

	which is used in \cite{sussmanLevelSetApproach1994}. Equation (\ref{eqn:distance_level_set}) can be discretized similarly as in (\ref{discr_level_set}). Note that $\phi$ at the grid points shifted to the PEC boundary will be set to 0 throughout the evolving of (\ref{eqn:distance_level_set}).

	With the signed distance function $\phi$ in mind, we can then compute the local normal direction $\bm{n}=\frac{\nabla \phi}{\| \nabla \phi \|}$ and tangential direction $\bm{t}$ which are used throughout the algorithm. To obtain $\bm{n}$, we use linear least squares module to fit locally a polynomial and then consider its gradient as the approximation of $\nabla \phi$. Its corresponding tangential direction is then obtained by rotating $\bm{n}$ clockwise by $\pi/ 2$.

	Now we introduce the approximation of ghost point values. Suppose we would like to extend $\bm{H}$ and construct its ghost field across the PEC boundary. First, we consider an orthogonal decomposition of $H^{gho}$, at the ghost point, along the local normal direction of the interface

	\begin{align}\label{eqn:ortho_mirr}
	\begin{split}
	\bm{H}^{gho} = (H_\perp)^{gho} \bm{n} + (H_\parallel)^{gho}  \bm{t}~.\\
	\end{split}
	\end{align}

	We then consider the first two terms of the Taylor expansions of $(H_\perp)^{gho}$ and $(H_\parallel)^{gho}$ along $\bm{n}$ around the point where the local normal line meets the PEC boundary,

	\begin{align}\label{taylor_mirr}
	\begin{split}
	(H_\perp)^{gho} = (H^{(0)}_\perp)^{gho} + ( \nabla H^{(0)}_\perp \cdot \bm{n} )^{gho} \phi \ldots\\
	(H_\parallel)^{gho} = (H^{(0)}_\parallel)^{gho} + ( \nabla  H^{(0)}_\parallel \cdot \bm{n} )^{gho} \phi \ldots\\
	\end{split}
	\end{align}

	where $\phi$ is the signed distance function between the ghost point and the PEC boundary, $\phi > 0$ inside the PEC object. It is now clear that to construct the ghost point values of the ghost fields, we need to obtain $ (H^{(0)}_\perp)^{gho}, ( \nabla H^{(0)}_\perp \cdot \bm{n} )^{gho}, \phi, (H^{(0)}_\parallel)^{gho} $ and $( \nabla  H^{(0)}_\parallel \cdot \bm{n} )^{gho}$. And these will be facilitated by the level set method.  

	In its original form in 2D, we consider level set function $\Phi$ and a closed curve defined as its zero level set $\Gamma = \{(x,y): \Phi(x,y) = 0\}$. To study the evolution of $\Gamma$, let $\Phi$ advect with field $\bm{v}$. Therefore, the level set equation \cite{OsherSethian88} is given by

	\begin{align}\label{eqn:level_set}
	\begin{split}
	\Phi_t + \bm{v} \cdot \nabla \Phi = 0~.\\
	\end{split}
	\end{align}

	In our case, we set $\Phi$ to be the components we like to extend across the interface following \cite{Peng99, Fedkiw99}. With the previous example, $\Phi$ will be replaced by $H_\perp$, $H_\parallel$, $\nabla H_\perp \cdot \bm{n}$ and $\nabla H_\parallel \cdot \bm{n}$, where the gradients are approximated by linear least squares as in (\ref{scheme:ls_cd}). Notice that the omission of the mirror superscripts. At every grid point outside the PEC object, decompose

	\begin{align}\label{ortho_inside}
	\begin{split}
	\bm{H} = (H_\perp) \bm{n} + (H_\parallel) \bm{t}~.\\
	\end{split}
	\end{align}

	These can be readily calculated at each given time in the algorithm progression. Next, we extend these components constantly along the local normal direction $\bm{n}$ (i.e. $\bm{v}$ is replaced by $\bm{n}$ in (\ref{eqn:level_set})) across the PEC boundary. Just like the discretization of Maxwell's Equations, we discretize the advection equation (\ref{eqn:level_set}) in the similar fashion

	\begin{align}\label{discr_level_set}
	\begin{split}
	\frac{\Phi^{n+1}_{i,j} - \tilde{\Phi}^n_{i,j} }{\Delta t} + \tilde{\bm{v}} \cdot \nabla \tilde{\Phi}^n_{i,j} = 0~,\\
	\end{split}
	\end{align}

	which gives us the update equation

	\begin{align}\label{eqn:update_eqn_level_set}
	\begin{split}
	\Phi^{n+1}_{i,j} = \tilde{\Phi}^n_{i,j} - \left( \tilde{\bm{v}} \cdot \nabla \tilde{\Phi}^n_{i,j} \right)  \Delta t~,\\
	\end{split}
	\end{align}

	where terms with tildes are approximated with the linear least squares method.

	\subsection{Algorithm}

	We break our process into the following 3 modules: Preparations, which includes the procedure to locally conform the grids to the boundary with Point-shifting and computation of distance function $\phi$, Ghost Point Extension, which extends $\bm{H}$ and $E$ values across the PEC boundary, and the BFECC module, which is BFECC method applied to the hybrid scheme covered in 2.2.

	\begin{figure}[H]
  \centering
  \begin{forest}
  forked edges,
  for tree={
    grow'=0,
    draw,
    align=c,
  },
  highlight/.style={
    thick,
    font=\sffamily\bfseries
  }
  [
    [Preparation
    	[Point-Shift
		]
		[Compute signed distance function $\phi$ \\with $\phi>0$ inside the PEC object
		]
		[Ghost point \\ extension module
			[Set $\bm{H}\cdot \bm{n}$ and $E$ to 0 at grid \\points shifted to the PEC boundary.
			]
			[Compute ghost point \\values of $\bm{H}$ and $E$.
			]
		]
    ]
    [BFECC at \\each time step
    	[Intermediate \\Forward Step
    		[Apply ghost point extension module
			]
			[Perform forward advection of \\Maxwell's Equations with scheme (\ref{scheme:ls_cd})
			]
		]
		[Intermediate \\Backward Step
			[Apply ghost point extension module
			]
			[Perform backward advection of Maxwell's \\Equations with the same scheme applied \\to time reversed Maxwell's Equations.
			]
		]
		[Forward Step
			[Compensate the solution with $\bm{e}$ (Sec.~\ref{Sec:BFECC})
			]
			[Apply ghost point extension module
			]
			[
			Perform forward advection with (\ref{scheme:ls_cd})]
		]	
    ]
  ]
\end{forest}
  \caption{Algorithm Organization}\label{fig:org}
\end{figure}
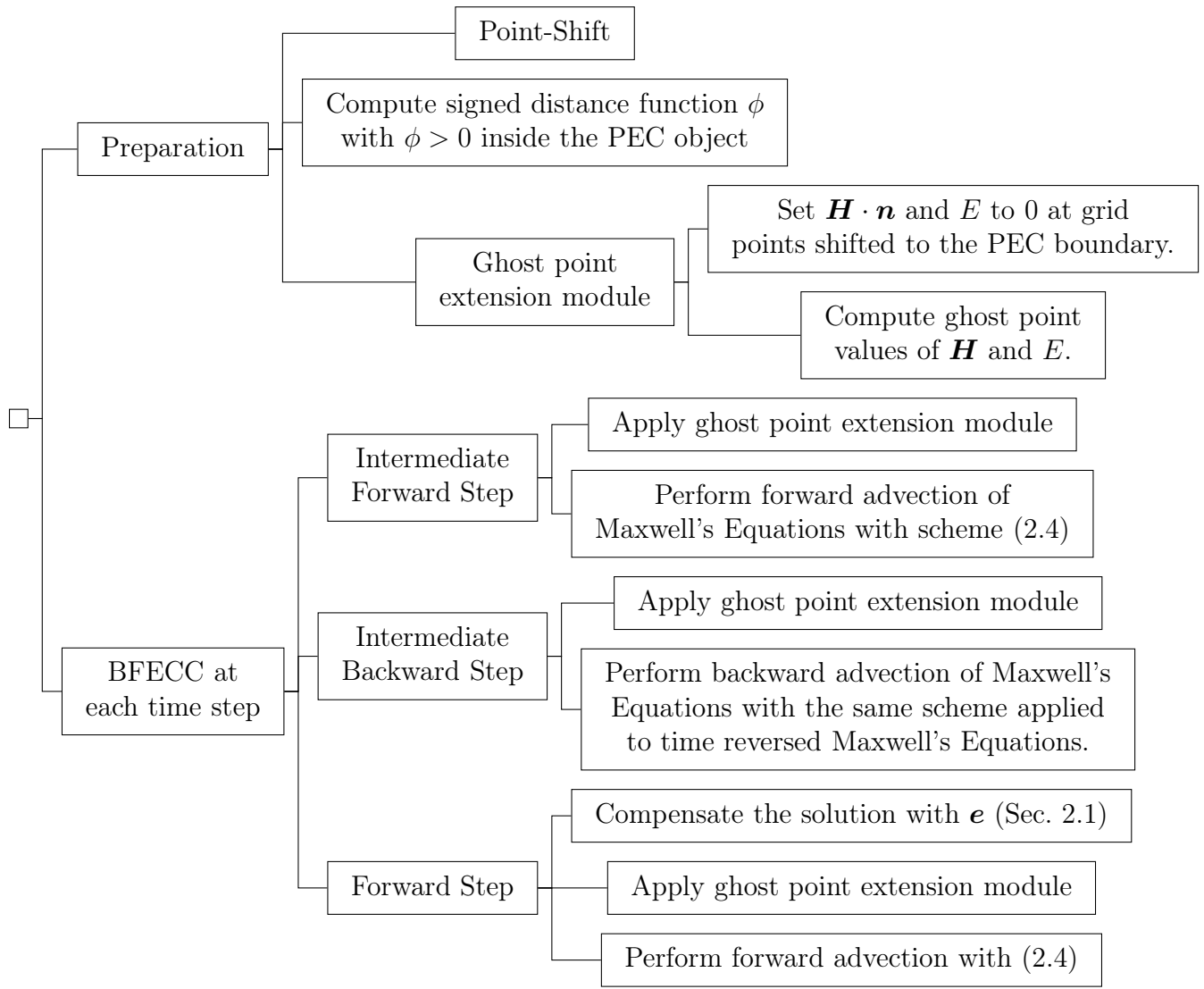
	
	Here we highlight the ghost point extension module for field $\bm{H}$. The extension of $E$ can be performed in a similar fashion.

\begin{algorithm}[H]\label{alg:point_shift}
    \SetKwInOut{Input}{Input}
    \SetKwInOut{Output}{Output}

    \Input{$\bm{H}$ values inside the computational domain at $t_n$, local normal direction $\bm{n}$, signed distance function $\phi$ to the PEC boundary}
    \Output{$\bm{H}$'s ghost point values inside PEC material at $t_n$}
    \begin{enumerate}

		\item Decompose $\bm{H}$ orthogonally into $H_\perp \bm{n}$ and $H_\parallel \bm{t}$.
		\item On the PEC boundary, set $H_\perp$ to be 0.
		\item In the computational domain, obtain $\nabla H_\perp \cdot \bm{n}$ and $\nabla H_\parallel \cdot \bm{n}$.
		\item Constantly extend $H^{(0)}_\perp, \nabla H^{(0)}_\perp \cdot \bm{n}, H^{(0)}_\parallel,  \nabla H^{(0)}_\parallel \cdot \bm{n}$ across the PEC boundary along the direction $\bm{n}$, where superscript $(0)$ indicates the values at the PEC boundary.
		\item With the above extensions, denote the corresponding values at ghost points  as $(H^{(0)}_\perp)^{gho}, ( \nabla H^{(0)}_\perp \cdot \bm{n} )^{gho}, (H^{(0)}_\parallel )^{gho}, ( \nabla H^{(0)}_\parallel \cdot \bm{n} )^{gho}$. Note that $(H^{(0)}_\perp)^{gho} = 0$ because $H^{(0)}_{\perp}$ at the PEC boundary is set to 0.
		\item Construct desired values at ghost points with linear parts of their Taylor expansions, that is

		$$(H_\perp)^{gho} = 0 + ( \nabla H^{(0)}_\perp \cdot \bm{n} )^{gho} \phi~, $$
		$$(H_\parallel)^{gho} = (H^{(0)}_\parallel)^{gho} - ( \nabla H^{(0)}_\parallel \cdot \bm{n} )^{gho} \phi~.$$

		\item Construct $\bm{H}$ at ghost points, that is 

		$\bm{H}^{gho} = (H_\perp)^{gho} \bm{n} + (H_\parallel)^{gho}  \bm{t}~.$

	\end{enumerate}
    \caption{Subroutine of ghost point extension for $\bm{H}$ at $t_n$}
\end{algorithm}

\medskip

\emph{\bf Remarks.} 

	\begin{enumerate}
                \item $\phi$ is the value of the signed distance function at the corresponding ghost point.
		\item The ghost point extension of the scalar field $E$ is identical to that of $H_\perp$.
		\item The extension of $H_\perp$ and $E$ is the so-called odd extension while the extension of $H_{\parallel}$ is an even extension.

		\item Since $\nabla H^{(0)}_\perp \cdot \bm{n}, \nabla H^{(0)}_\parallel \cdot \bm{n}$ are not easy to compute at the PEC boundary, we actually compute $\nabla H_\perp \cdot \bm{n}, \nabla H_\parallel \cdot \bm{n}$ at grid points adjacent to the PEC boundary and then constantly extend them along local normal direction $\bm{n}$ across the PEC boundary until well passing the band formed by ghost points. This will not hurt the order of accuracy because we only need the gradients to be first order accurate.

		\item Only one layer of ghost points beyond the PEC boundary is needed for the underlying scheme of BFECC to run.

		\item In order to ensure a constant extension of boundary values along local normal direction $\bm{n}$, the values of the components to be extended ($H_\perp$ and $H_\parallel$ at grid points where $\phi \le 0$; $\nabla H_\perp \cdot \bm{n}$ and $\nabla H_\parallel \cdot \bm{n}$ at grid points where $\phi < 0$) must remain unchanged and not be updated by scheme (\ref{discr_level_set}).
	\end{enumerate}

\section{Numerical experiments}\label{sec:num_experiments}
\subsection{Cylindrical Perfect Electric Conductor with a  Circular Cross Section}

We first consider a circular cylindrical perfect electric conductor, whose axis extends along $z$-direction, subjected to a $z$-polarized incident plane wave traversing along the positive direction of $x$-axis. The incident wave is given by $E_z^{i} = sin( \omega (x-t)), H_x^{i} = 0$ and $H_y^{i} = - sin( \omega (x - t))$, where angular frequency $\omega = 2\pi/0.6$. The right circular cylindrical PEC is assumed to be of negligible thickness. The cross section of the cylinder by $x$-$y$ axis is given by a circle, centered at $(5.0, 5.0)$ with radius $r = 2.0$. The medium that encloses the cylinder has $\epsilon = 1$ and $ \mu = 1$  as its permittivity and permeability. The computational domain is given by $[0, 10] \times [0, 10]$. An illustration is given by the following figure. 

\begin{figure}[H]
\centering
\includegraphics[width=1\linewidth]{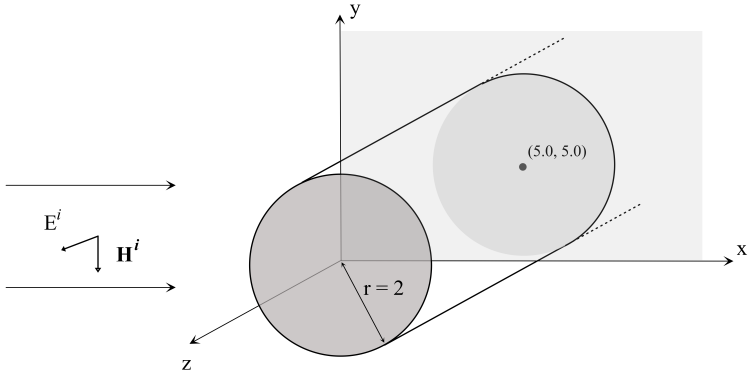}
\caption{Experiment setup}
\label{fig:Experiment setup}
\end{figure}

We solve the system with our method from $T=0$ to $T=2$, with $\Delta t/ \Delta x = 1$ on grids of sizes $100 \times 100$, $200 \times 200$, $400 \times 400$, $800 \times 800$, $1600 \times 1600$ and $3200 \times 3200$. Further, we use the solution obtained on the $3200 \times 3200$ as our exact solution to interpolate convergence trends. 

\begin{figure}[H]
\centering
\parbox{5cm}{
\includegraphics[width=5cm]{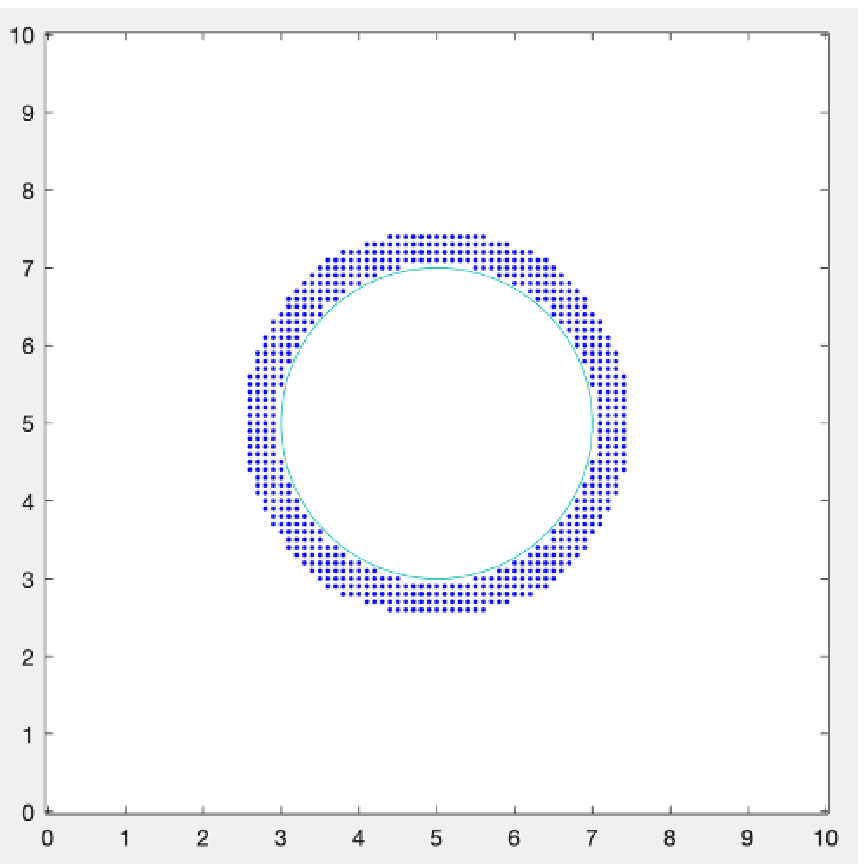}

\label{fig:Circular_error}}
\qquad\qquad
\begin{minipage}{5cm}
\includegraphics[width=5cm]{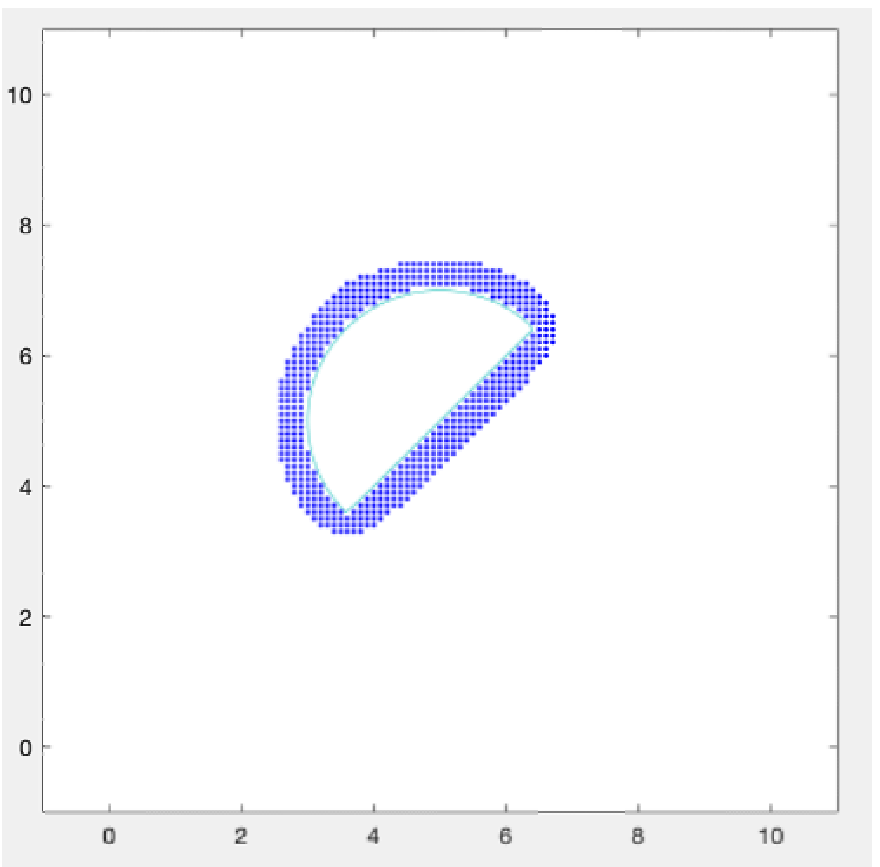}

\label{fig:Halfmoon_error}
\end{minipage}
\caption{Sampling points for grid refinement analysis: within 10 $\Delta x$ right outside the circular/half-moon PEC boundary on the coarsest mesh}\label{fg:error_sampling}
\end{figure}

For the grid refinement analysis, we opted to use the grid points within $10 \Delta x$ from the PEC boundary at the coarsest mesh. This is because at $T=1$, the wave has propagated over a distance of $10 \Delta x$ at the coarsest mesh. This choice is also used in the analysis of the half-moon shaped PEC boundary example. Figure \ref{fg:error_sampling} provides a visual representation of such grid points. The order of accuracy result is given by Table \ref{tb:1d_uniform_error_1}, in which $\Delta t$ is chosen with $\Delta t / \Delta x = 1$ with $\Delta x$ being the mesh size before applying the point-shift procedure.

\begin{table}[H]
\centering
\begin{tabular}{|c|c|c|c|c|c|c|c|c|} \hline
 $\Delta t / \Delta x = 1$ & \multicolumn{4}{|c|}{ Measured in $l_1$ norms } \\ \hline
 & \multicolumn{2}{|c|}{ Measured in $E_z$ } & \multicolumn{2}{|c|}{ Measured in $B_x$ } \\ \hline
 Grid & Error & Order & Error & Order \\ \hline
$100$ & $6.69 \times 10^{-1}$ & -- & $3.20 \times 10^{-1}$ & -- \\ \hline
$200$ & $3.41 \times 10^{-1}$ & 0.96 & $1.59 \times 10^{-1}$ & 1.01 \\ \hline
$400$ & $1.06 \times 10^{-1}$  & 1.70 & $5.05 \times 10^{-2}$ & 1.66 \\ \hline
$800$ & $2.71 \times 10^{-2}$  & 1.96 & $1.39 \times 10^{-2}$ & 1.86 \\ \hline
\end{tabular}
\caption{Circular cross-section: Order of accuracy for BFECC based on the central difference scheme at $T = 1$}\label{tb:1d_uniform_error_1}
\end{table}

\begin{figure}[H]
\centering
\includegraphics[width=1\linewidth]{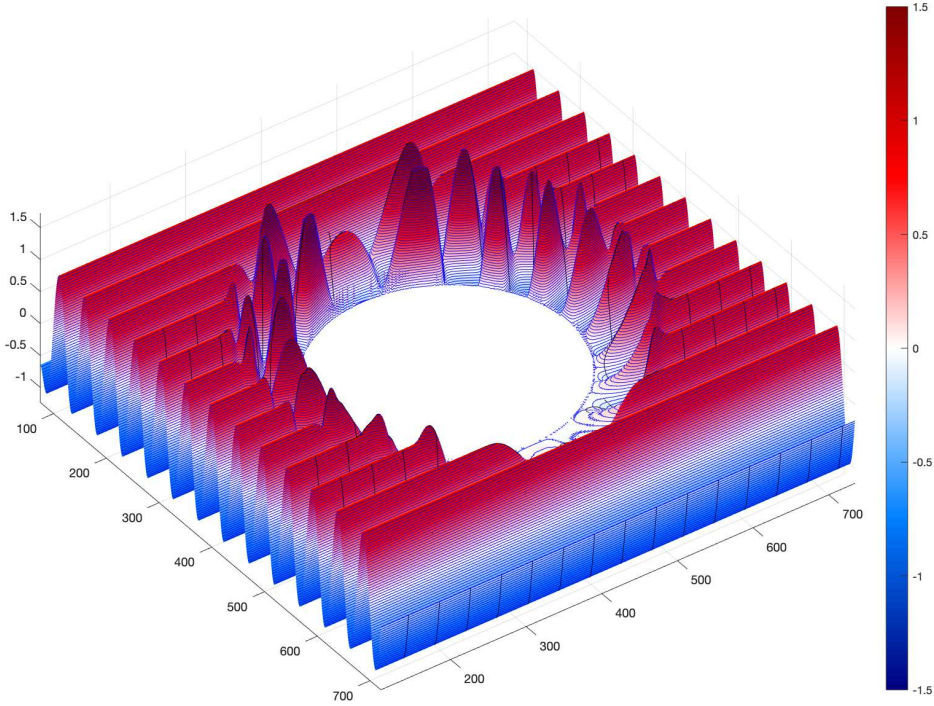}
\caption{Graph of $E_z$, with circular cross-section}
\label{fig:graph of E_z}
\end{figure}

\begin{figure}[H]
\centering
\includegraphics[width=0.8\linewidth]{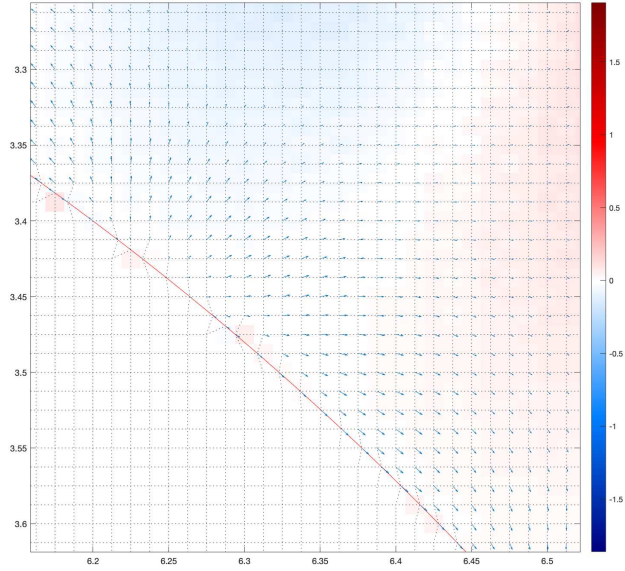}
\caption{Northeast Corner of the PEC Cross-section, with field $\bm{H}$ and $E_z$}
\label{fig:Northeast corner}
\end{figure}

Figure \ref{fig:graph of E_z} is the graph of $E_z$ from an angle. We can observe how the wave "bounces" off against the cylindrical PEC object. Figure \ref{fig:Northeast corner} provides a close-up view of the northeast corner of the PEC cross-section. The blue vectors represent field $\bm{H}$ whereas the blue-red color reflects magnitude of $E_z$. One can also observe how field $\bm{H}$ is behaving according to the PEC condition near the PEC boundary.

Remarks:

1. Note that when updating gradients on the PEC boundary, for example with $\nabla E$, we consider the interior points within the computational domain to extend values to ghost point locations inside the PEC medium. This is because if we instead use grid points on the boundary, we will only obtain cumulatively $O(1)$ accuracy, which will lower the overall accuracy of the scheme.\\

2. To obtain the distance function $\phi$ to the boundary, we run a similar scheme over fictitious time $t$ until $\| \nabla \phi \|$ stabilizes and converges to $1$. In this and later numerical experiments, we have computed the distance function over the whole computational domain. However, to speed up the process, one only needs to determine $\phi$ for grid points within a thin band with the PEC boundary at the center for the extension module to work. The scheme for computing $\phi$ in our experiment uses $\Delta t / \Delta x = 0.2$.\\

3. For the hybrid scheme with BFECC applied, we iterate the module until the wave has traversed beyond $10 \Delta x$ away from the boundary. The grids points within $10 \Delta x$ are later used to determine the accuracy of the scheme.

\subsection{Cylindrical Perfect Electric Conductor with a half-moon cross section}

\begin{figure}[H]
\centering
\includegraphics[width=0.8\linewidth]{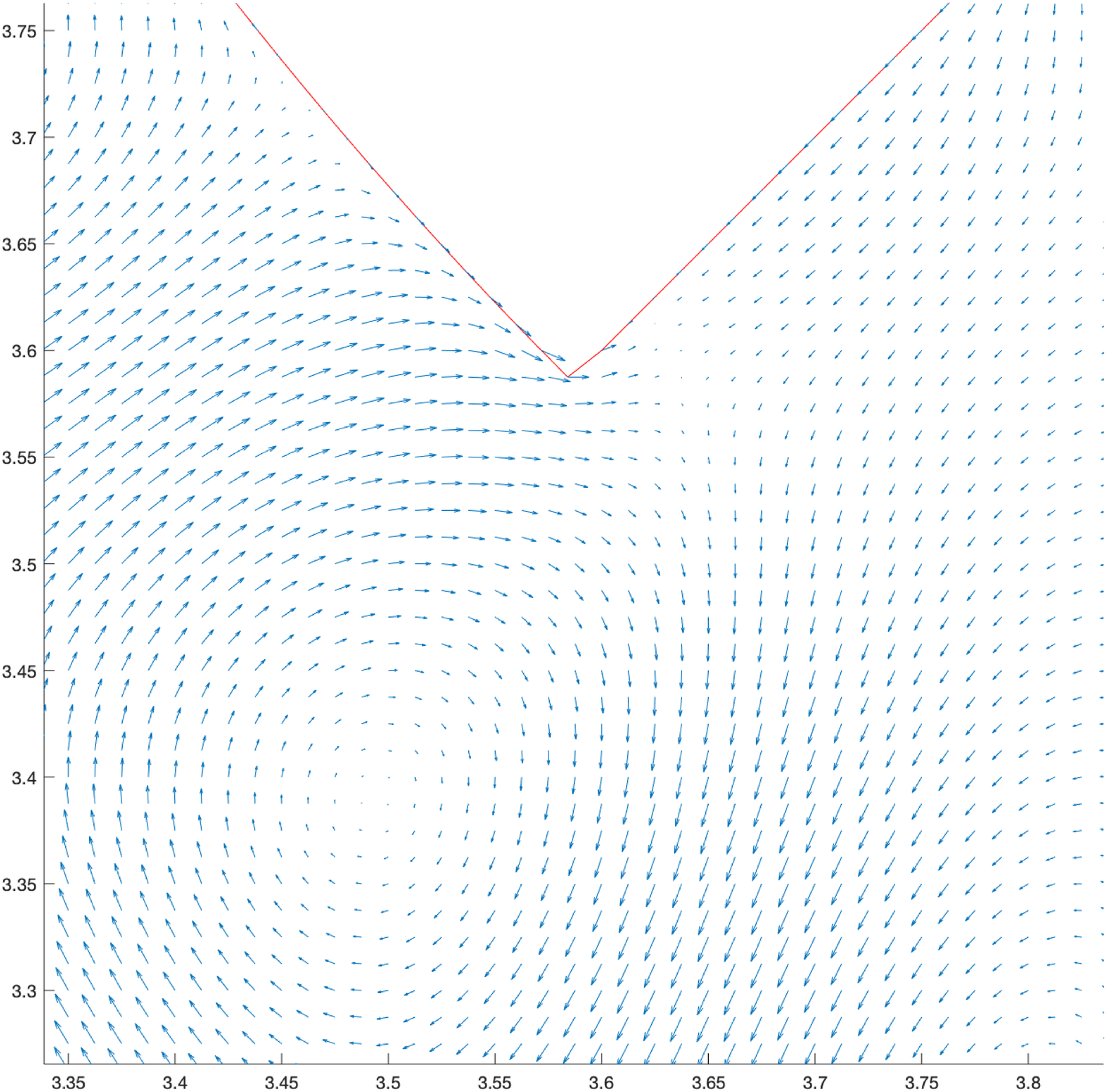}
\caption{No artifacts are observed near the bottom corner of the half-moon PEC boundary}
\label{fig:Halfmoon_error_closeup}
\end{figure}

The reason why we chose to test this scenario is because we would like to see how our method would behave around sharp corners with point-shifted grids. and the sharp corners may affect order of accuracy. Things to watch out for in here would be computational artifacts such as sudden jumps or turbulence. Figure \ref{fig:Halfmoon_error_closeup} highlights a close-up of the bottom corner of the half-moon cross-section, where tiny blue vectors are field $\bm{H}$. Fields near the top corner also behave in a similar way.

We do not observe any artifacts in terms of field $\bm{H}$ near the sharp corners of the half-moon PEC boundary. A 3D level sets for $E_z$ is given in Figure \ref{fig:graph of E_z, half-moon cross-section}. No artifacts are observed in terms of $E_z$, either.

\begin{figure}[H]
\centering
\includegraphics[width=1\linewidth]{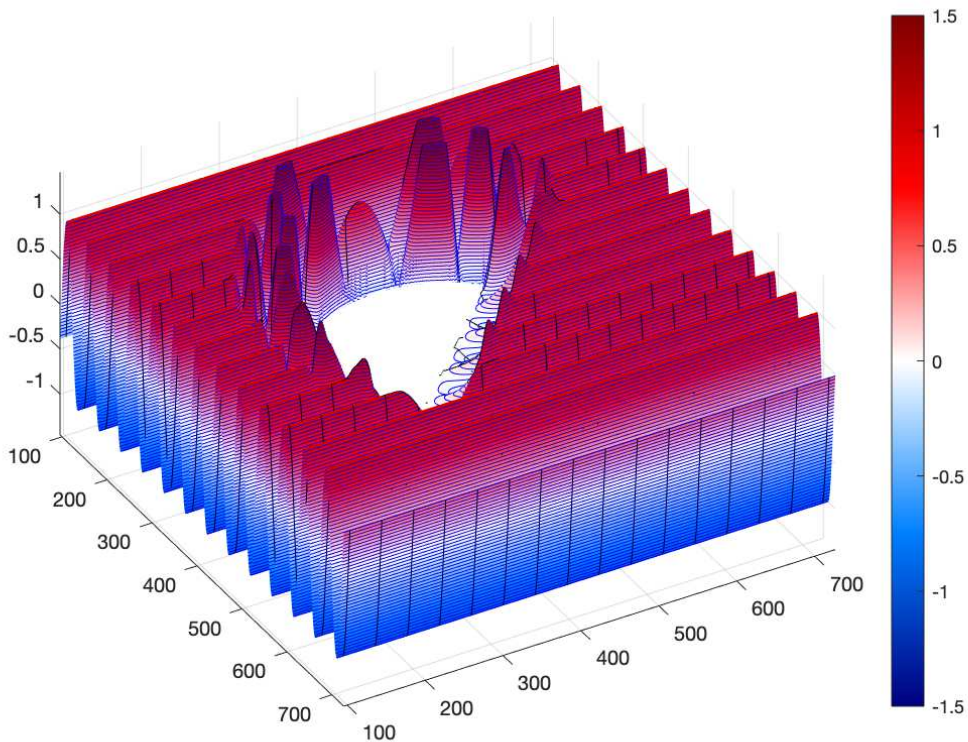}
\caption{Graph of $E_z$, with a half-moon cross-section}
\label{fig:graph of E_z, half-moon cross-section}
\end{figure}

As for accuracy, Table \ref{tb:half_moon_1d_uniform_error_1} shows that the same converging trend is still present, not very different from the cylindrical case.

\begin{table}[H]
\centering
\begin{tabular}{|c|c|c|c|c|} \hline
 $\Delta t / \Delta x = 1$ & \multicolumn{4}{|c|}{ Measured in $l_1$ norms } \\ \hline
 & \multicolumn{2}{|c|}{ Measured in $E_z$ } & \multicolumn{2}{|c|}{ Measured in $B_x$ } \\ \hline
 Grid & Error & Order & Error & Order \\ \hline
$100$ & $6.20 \times 10^{-1}$ & -- & $3.01 \times 10^{-1}$ & --\\ \hline
$200$ & $3.24 \times 10^{-1}$ & 0.93 & $1.47 \times 10^{-1}$ & 1.03 \\ \hline
$400$ & $1.00 \times 10^{-2}$  & 1.70 & $4.72 \times 10^{-2}$ & 1.64 \\ \hline
$800$ & $2.54 \times 10^{-2}$  & 1.98 & $1.31 \times 10^{-2}$ & 1.85 \\ \hline
\end{tabular}
\caption{Half-Moon Cross-section: Order of accuracy for BFECC based on the central difference scheme at $T = 1$}\label{tb:half_moon_1d_uniform_error_1}
\end{table}

\subsection{Cylindrical Perfect Electric Conductor with a circular cross-section where $\Delta t / \Delta x = 1.4$}

The reason why we tested the scenario where $\Delta t / \Delta x = 1.4$ is because in \cite{wangBackForthError2019}, it is proved that in a two dimensional case, BFECC can provide second order accuracy applied to a central difference scheme over an orthogonal grid when $\Delta t \le \sqrt{3}/ (( 1/\Delta x )^2 + ( 1/\Delta y )^2 )$. Translating this into our case, the CFL condition requires that $\Delta t / \Delta x $ is at most $\sqrt{3/2}$. Although the underlying grid is not orthogonal because of the point-shifting, it is reasonable to believe that at $\Delta t / \Delta x = 1.4$ we are really creating a strenuous condition that would potentially strongly undermine the performance of the scheme. However, Table \ref{tb:1d_uniform_error_1_1.4} shows that the converging trend is still strong.

\begin{table}[H]
\centering
\begin{tabular}{|c|c|c|c|c|} \hline
 $\Delta t / \Delta x = 1.4$ & \multicolumn{4}{|c|}{ Measured in $l_1$ norms } \\ \hline
 & \multicolumn{2}{|c|}{ Measured in $E_z$ } & \multicolumn{2}{|c|}{ Measured in $B_x$ } \\ \hline
 Grid & Error & Order & Error & Order \\ \hline
$200$ & $5.86 \times 10^{-1}$ & -- & $2.57 \times 10^{-1}$ & -- \\ \hline
$400$ & $2.64 \times 10^{-1}$  & 1.15 & $1.20 \times 10^{-1}$ & 1.09 \\ \hline
$800$ & $6.84 \times 10^{-2}$  & 1.95 & $3.24 \times 10^{-2}$ & 1.89 \\ \hline
\end{tabular}
\caption{Circular Cross-section at $T = 1$, $\Delta t / \Delta x = 1.4$}\label{tb:1d_uniform_error_1_1.4}
\end{table}

As can be observed, the convergence rate is not as fast as in the case when $\Delta t / \Delta x = 1$. However, it is still fairly impressive as we are over the theoretical upper limit of the CFL condition. The fact that we were able to obtain relatively good results under such strenuous condition goes to show that the novel approach is indeed fairly robust.

\section{Conclusion}

Our method is shown to be effective in relation to the computational cost, strain and the ease of implementation. In fact, all experiments have been performed on a personal computer with minimal optimization. It does not require a careful analysis of the PEC object's geometry and can be easily extended to an arbitrary shape with an interior with minimal modification. BFECC applied to the underlying scheme improves the order of accuracy both in space and time. Its robustness has been highlighted even with a strenuous CFL condition.

\bibliographystyle{amsplain}

\bibliography{BFECC}

\end{document}